\documentclass{my_cmmse2010}

\usepackage{amssymb}
\usepackage{amsmath,amsfonts,amsthm}
\usepackage{graphicx}
\usepackage{url}
\usepackage{setspace}

\newtheorem{thm}{Proposition}

\begin{document}

\markboth
  {Control of dengue disease: a case study in Cape Verde} 
  {Rodrigues, Monteiro, Torres \& Zinober}


\title{Control of dengue disease:\\
a case study in Cape Verde}

\author{Helena Sofia Rodrigues}{sofiarodrigues@esce.ipvc.pt}{1}
\author{M. Teresa T. Monteiro}{tm@dps.uminho.pt}{2}
\author{Delfim F. M. Torres}{delfim@ua.pt}{3}
\author{Alan Zinober}{a.zinober@sheffield.ac.uk}{4}

\affiliation{1}{School of Business Studies}{Viana do Castelo Polytechnic Institute, Portugal}
\affiliation{2}{Department of Production and Systems}{University of Minho, Portugal}
\affiliation{3}{Department of Mathematics}{University of Aveiro, Portugal}
\affiliation{4}{Department of Applied Mathematics}{University of Sheffield, UK}


\begin{abstract}
A model for the transmission of dengue disease is presented. 
It consists of eight mutually-exclusive compartments representing 
the human and vector dynamics. It also includes a control parameter 
(adulticide spray) in order to combat the mosquito. The model presents 
three possible equilibria: two disease-free equilibria (DFE) --- where humans, 
with or without mosquitoes, live without the disease --- and another 
endemic equilibrium (EE). In the literature it has been proved that a DFE  
is locally asymptotically stable, whenever a certain epidemiological threshold, 
known as the \emph{basic reproduction number}, is less than one. 
We show that if a minimum level of insecticide is applied, 
then it is possible to maintain the basic reproduction number below unity. 
A case study, using data of the outbreak 
that occured in 2009 in Cape Verde, is presented.

\keywords dengue, basic reproduction number, stability, Cape Verde, control.
\msccodes 92B05, 93C95, 93D20.
\end{abstract}


\section{Introduction}

Dengue is a mosquito-borne infection that has become a major international 
public health concern. According to the World Health Organization, 
50 to 100 million dengue infections occur yearly, including 500000 
Dengue Haemorrhagic Fever cases and 22000 deaths, 
mostly among children \cite{WHO}. Dengue is found in tropical 
and sub-tropical regions around the world, 
predominantly in urban and semi-urban areas.

There are two forms of dengue: Dengue Fever and Dengue Haemorrhagic Fever. 
The first one is characterized by a sudden fever without respiratory symptoms, 
accompanied by intense headaches and lasts between three and seven days. 
The second one has the previous symptoms but also nausea, vomiting, 
fainting due to low blood pressure and can lead to death 
in two or three days \cite{Derouich2003}.

The spread of dengue is attributed to expanding geographic distribution 
of the four dengue viruses and their mosquito vectors, the most important 
of which is the predominantly urban species \emph{Aedes aegypti}. 
The life cycle of a mosquito presents four distinct stages: 
egg, larva, pupa and adult. In the case of \emph{Aedes aegypti} 
the first three stages take place in or near water while air 
is the medium for the adult stage \cite{Otero2008}. The adult stage 
of the mosquito is considered to last an average of eleven days 
in the urban environment. Dengue is spread only by adult females, 
that require a blood meal for the development of eggs; male mosquitoes 
feed on nectar and other sources of sugar. In this process the female 
acquire the virus while feeding on the blood of an infected person. 
After virus incubation for eight to ten days, an infected mosquito 
is capable, during probing and blood feeding, 
of transmitting the virus for the rest of its life.

The organization of this paper is as follows. 
A mathematical model of the interaction between human and mosquito populations 
is presented in Section~\ref{sec:2}. Section~\ref{sec:3} is concerned with 
the equilibria of the epidemiological model and their stability. 
In Section~\ref{sec:4} the results obtained in the previous section 
are applied to a study case. Finally, some concluding 
notes are given in Section~\ref{sec:5}.


\section{The mathematical model}
\label{sec:2}

Considering the work of \cite{Rodrigues2009}, the relationship between humans 
and mosquitoes are now rather complex, taking into account the model presented 
in \cite{Dumont2008}. The novelty in this paper is the presence 
of the control parameter related to adult mosquito spray.

The notation used in our mathematical model 
includes four epidemiological states for humans:

\begin{quote}
\begin{tabular}{ll}
$S_h(t)$ & susceptible (individuals who can contract
the disease)\\
$E_h(t)$ & exposed (individuals who have been infected by the parasite\\
& but are not yet able to transmit to others)\\
$I_h(t)$ & infected (individuals capable of transmitting the disease to others)\\
$R_h(t)$ & resistant (individuals who have acquired immunity)\\
\end{tabular}
\end{quote}

It is assumed that the total human population $(N_h)$ is constant,
so, $N_h=S_h+E_h+I_h+R_h$.
There are also other four state variables related to the female
mosquitoes (the male mosquitoes are not considered in this study
because they do not bite humans and consequently they do not
influence the dynamics of the disease):

\begin{quote}
\begin{tabular}{ll}
$A_m(t)$& aquatic phase (that includes the egg, larva and pupa stages)\\
$S_m(t)$& susceptible (mosquitoes that are able to contract the disease)\\
$E_m(t)$& exposed (mosquitoes that are infected but are not yet able \\
& to transmit to humans)\\
$I_m(t)$& infected (mosquitoes capable of transmitting the disease to humans)\\
\end{tabular}
\end{quote}

In order to analyze the effects of campaigns to combat the mosquito, 
there is also a control variable:

\begin{quote}
\begin{tabular}{ll}
$c(t)$& level of insecticide campaigns\\
\end{tabular}
\end{quote}

Some assumptions are made in this model:
\begin{itemize}
\item the total human population ($N_h$) is constant, 
which means that we do not consider births and deaths;
\item there is no immigration of infected individuals 
to the human population;
\item the population is homogeneous, which means that every 
individual of a compartment is homogenously mixed with the other individuals;
\item the coefficient of transmission of the disease is fixed and do not vary seasonally;
\item both human and mosquitoes are assumed to be born susceptible; 
there is no natural protection;
\item for the mosquito there is no resistant phase, 
due to its short lifetime.
\end{itemize}

The parameters used in our model are:

\begin{quote}
\begin{tabular}{ll}
$N_h$ & total population \\
$B$ & average daily biting (per day)\\
$\beta_{mh}$ & transmission probability from $I_m$ (per bite) \\
$\beta_{hm}$ & transmission probability from $I_h$ (per bite) \\
$1/\mu_{h}$ & average lifespan of humans (in days) \\
$1/\eta_{h}$ & mean viremic period (in days)\\
$1/\mu_{m}$ & average lifespan of adult mosquitoes (in days) \\
$\mu_{b}$ & number of eggs at each deposit per capita (per day) \\
$\mu_{A}$ & natural mortality of larvae (per day) \\
$\eta_{A}$ & maturation rate from larvae to adult (per day) \\
$1/\eta_{m}$ & extrinsic incubation period (in days)  \\
$1/\nu_{h}$ & intrinsic incubation period (in days) \\
$m$ & female mosquitoes per human \\
$k$ & number of larvae per human \\
$K$ & maximal capacity of larvae\\
\end{tabular}
\end{quote}

The Dengue epidemic can be modelled by the following nonlinear
time-varying state equations:

Human Population
\begin{equation}\label{odehuman}
\begin{tabular}{l}
$
\left\{
\begin{array}{l}
\frac{dS_h}{dt}(t) = \mu_h N_h - (B\beta_{mh}\frac{I_m}{N_h}+\mu_h) S_h\\
\frac{dE_h}{dt}(t) = B\beta_{mh}\frac{Im}{N_h}S_h - (\nu_h + \mu_h )E_h\\
\frac{dI_h}{dt}(t) = \nu_h E_h -(\eta_h  +\mu_h) I_h\\
\frac{dR_h}{dt}(t) = \eta_h I_h - \mu_h R_h\\
\end{array}
\right. $\\
\end{tabular}
\end{equation}

and vector population

\begin{equation}\label{odevector}
\begin{tabular}{l}
$
\left\{
\begin{array}{l}
\frac{dA_m}{dt}(t) = \mu_b (1-\frac{A_m}{K})(S_m+E_m+I_m)-(\eta_A+\mu_A) A_m\\
\frac{dS_m}{dt}(t) = -(B \beta_{hm}\frac{I_h}{N_h}+\mu_m) S_m+\eta_A A_m-c S_m\\
\frac{dE_m}{dt}(t) = B \beta_{hm}\frac{I_h}{N_h}S_m-(\mu_m + \eta_m) E_m-c E_m\\
\frac{dI_m}{dt}(t) = \eta_m E_m -\mu_m I_m - c I_m\\
\end{array}
\right. $\\
\end{tabular}
\end{equation}
with the initial conditions
\begin{equation}
\label{initial}
\begin{tabular}{llll}
$S_h(0)=S_{h0},$ & $E_h(0)=E_{h0},$ & $I_h(0)=I_{h0},$ &
$R_h(0)=R_{h0},$ \\
$A_m(0)=A_{m0},$ & $S_{m}(0)=S_{m0},$ &
$E_m(0)=E_{m0},$ & $I_m(0)=I_{m0}.$
\end{tabular}
\end{equation}

Notice that the equation related to the aquatic phase does not have
the control variable $c$, because the adulticide does not produce
effects in this stage of the life of the mosquito.


\section{Equilibrium points and Stability}
\label{sec:3}

Let the set
\begin{center}
\begin{tabular}{l}
\small
$\Omega=\{(S_h,E_h,I_h,A_m,S_m,E_m,I_m)\in \mathbb{R}^{7}_{+}: 
S_h+E_h+I_h\leq N_h, A_m\leq k N_h, S_m+E_m+I_m\leq m N_h \}$
\end{tabular}
\end{center}

\noindent be the region of biological interest, that is positively invariant 
under the flow induced by the differential system (\ref{odehuman})--(\ref{odevector}).

\begin{thm}
Let $\Omega$ be defined as above. Consider also
\begin{center}
\begin{tabular}{l}
$\mathcal{M}=-\left(c (\eta_A + \mu_A) + \mu_A \mu_m + \eta_A (-\mu_b + \mu_m)\right)$.
\end{tabular}
\end{center}

The system (\ref{odehuman})--(\ref{odevector}) admits, at most, three equilibrium points:
\begin{itemize}
\item if $\mathcal{M}\leq 0$, there is a Disease-Free Equilibrium (DFE), 
called Trivial Equilibrium, $E_{1}^{*}=\left(N_h,0,0,0,0,0,0\right)$;
\item if $\mathcal{M}> 0$, there is a Biologically 
Realistic Disease-Free Equilibrium (BRDFE),

$E_{2}^{*}=\left(N_h,0,0,\frac{k Nh \mathcal{M}}{\eta_A\mu_b},
\frac{k Nh \mathcal{M}}{\mu_b \mu_m},0,0\right)$

or an Endemic Equilibrium (EE), 
$E_{3}^{*}=\left(S_h^*,E_h^*,I_h^*,A_m^*,S_m^*,E_m^*,I_m^*\right)$.

\end{itemize}
\end{thm}

It is necessary to determine the \emph{basic reproduction number} of the disease, 
$\mathcal{R}_0$. This number is very important from the epidemiologistic 
point of view. It represents the expected number of secondary cases produced 
in a completed susceptible population, by a typical infected individual 
during its entire period of infectiousness \cite{Hethcote2000}.
Following \cite{Driessche2002}, we prove:

\begin{thm}
If $\mathcal{M}>0$, then the basic reproduction 
number associated to (\ref{odehuman})--(\ref{odevector}) 
is $\mathcal{R}_{0}^2=\frac{B^2 k \beta_{hm} \beta_{mh} 
\eta_m \nu_h\mathcal{M} }{\mu_b (\eta_h + \ \mu_h) (c 
+ \mu_m)^2 (c + \eta_m + \mu_m) (\mu_h + \nu_h)}$.

\medskip
BRDFE is locally asymptotically stable if $\mathcal{R}_{0}<1$ 
and unstable if $\mathcal{R}_{0}>1$.
\end{thm}

From a biological point of view, it is desirable that humans 
and mosquitoes coexist without the disease reaching a level 
of endemicity. We claim that proper use of the control $c$ 
can result in the basic reproduction number remaining 
below unity and, therefore, making BRDFE stable.

In order to make effective use of achievable insecticide control, 
and simultaneously to explain more easily to the competent authorities 
its effectiveness, we assume that $c$ is constant.

We want to find $c$ such that $\mathcal{R}_{0}<1$.


\section{Dengue in Cape Verde}
\label{sec:4}

The simulations were carried out using the following values: 
$N_h=480000$, $B=1$, $\beta_{mh}=0.375$,
$\beta_{hm}=0.375$, $\mu_{h}=1/(71*365)$, $\eta_{h}=1/3$,
$\mu_{m}=1/11$, $\mu_{b}=6$, $\mu_{A}=1/4$, $\eta_A=0.08$,
$\eta_m=1/11$, $\nu_h=1/4$, , $m=6$,
$k=3$, $K=k*N_h$.
The initial conditions for the problem were: $S_{h0}=m*N_h$,
$E_{h0}=216$, $I_{h0}=434$, $R_{h0}=0$, $A_{m0}=k*N_h$,
$S_{m0}=m*N_h$, $E_{m0}=0$, $I_{m0}=0$. The final time was
$t_f=84$ days. The values related to humans describes the
reality of an infected period
in Cape Verde \cite{CDC2010}. However, since it was the
first outbreak that happened in the archipelago it was not possible
to collect any data for the mosquito. Thus, for the \emph{aedes
Aegypti} we have selected information from Brazil where dengue is
already a reality long known \cite{Thome2010, Yang2009}.

\begin{thm}
Let us consider the parameters listed above and consider 
$c$ as a constant. Then $\mathcal{R}_{0}<1$ if and only if $c>0.0837$.
\end{thm}

For our computations let us consider $c=0.084$. The results indicate that
use of the control $c$ is crucial to prevent that an outbreak could transform 
an epidemiological episode to an endemic disease. The computational 
experiences were carried out using Scilab \cite{Campbell2006}.

\begin{figure}[ptbh]
\centering
\begin{minipage}[t]{0.48\linewidth}
\centering
\includegraphics[scale=0.5]{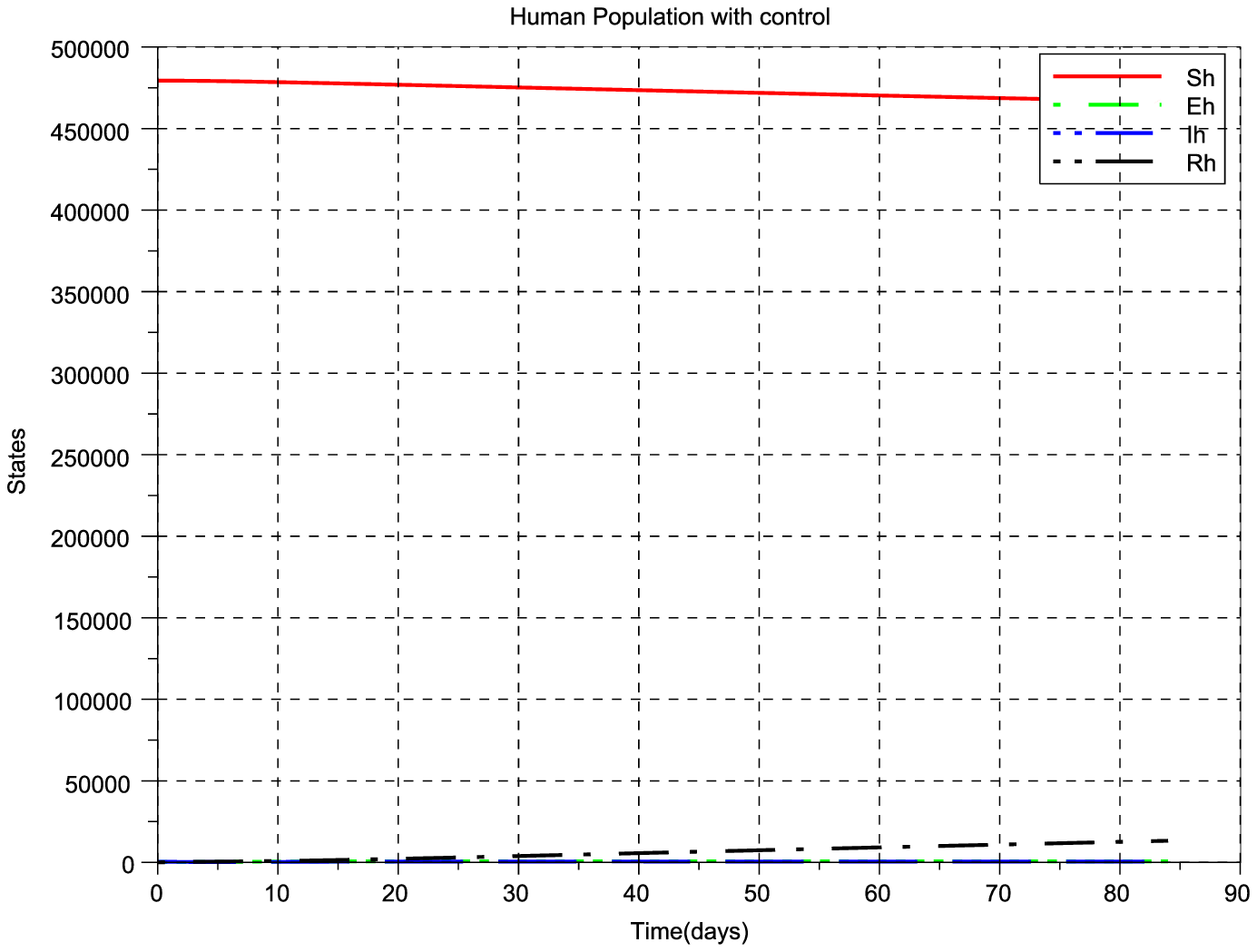}
{\caption{\label{human_control}  Human compartments using control}}
\end{minipage}\hspace*{\fill}
\begin{minipage}[t]{0.48\linewidth}
\centering
\includegraphics[scale=0.5]{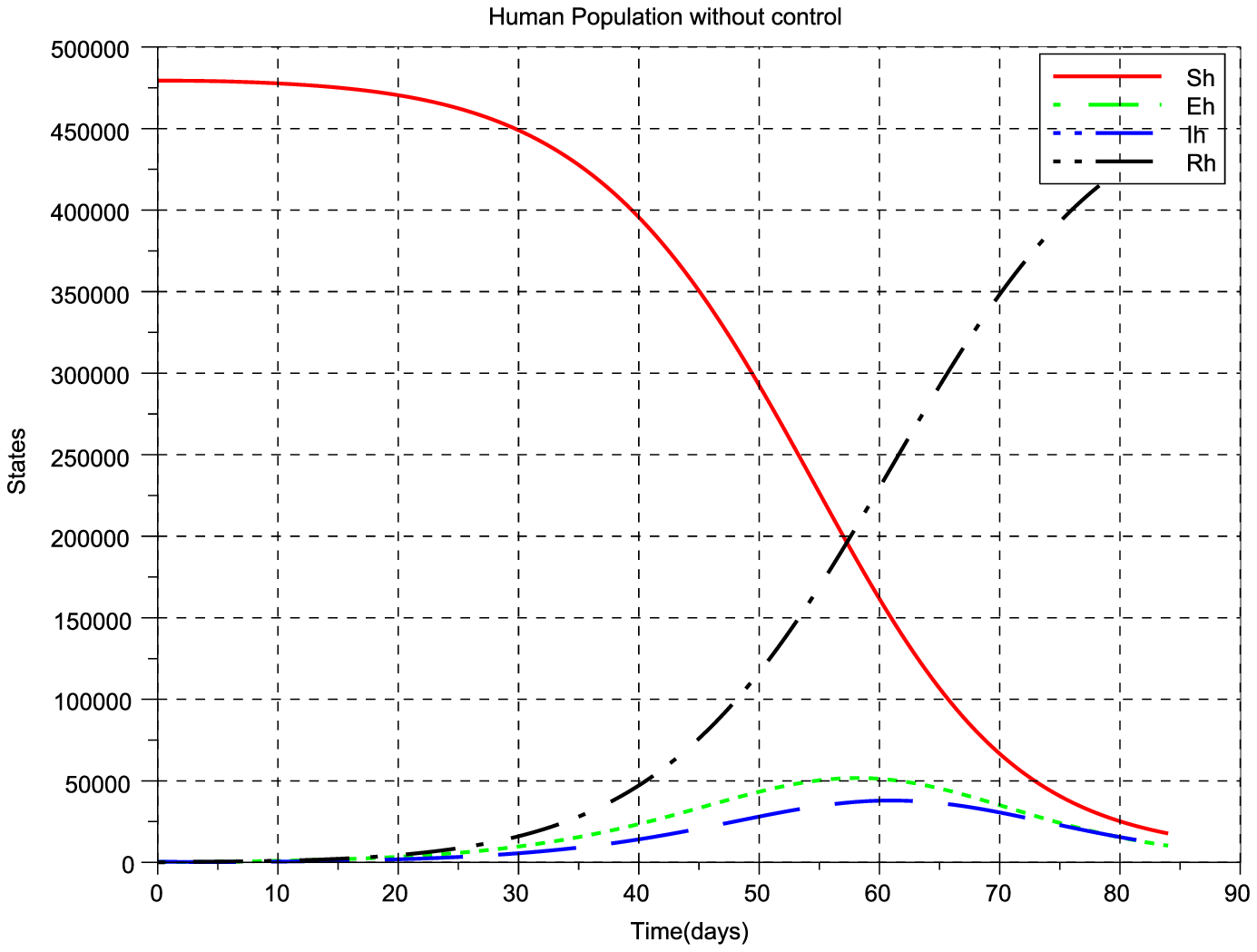}
{\caption{\label{human_nocontrol} \small Human compartments with no control}}
\end{minipage}
\end{figure}

Figures \ref{human_control} and \ref{human_nocontrol} show the curves related 
to human population, with and without control, respectively. The number 
of infected persons, even with small control, 
is much less than without any spray campaign.


\begin{figure}[ptbh]
\centering
\begin{minipage}[t]{0.48\linewidth}
\centering
\includegraphics[scale=0.50]{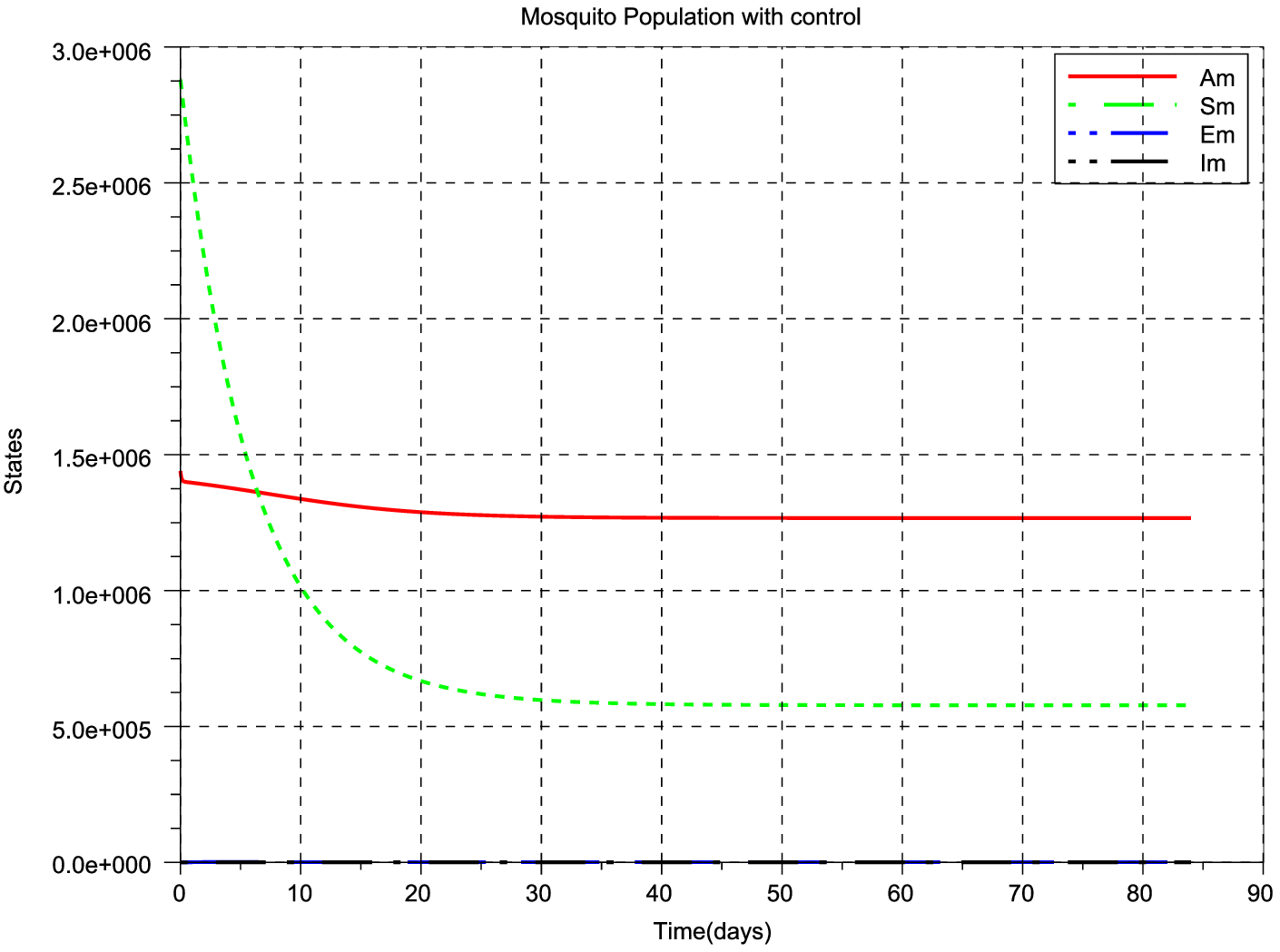}
{\caption{\label{mosquito_control}  Mosquito compartments using control}}
\end{minipage}\hspace*{\fill}
\begin{minipage}[t]{0.48\linewidth}
\centering
\includegraphics[scale=0.50]{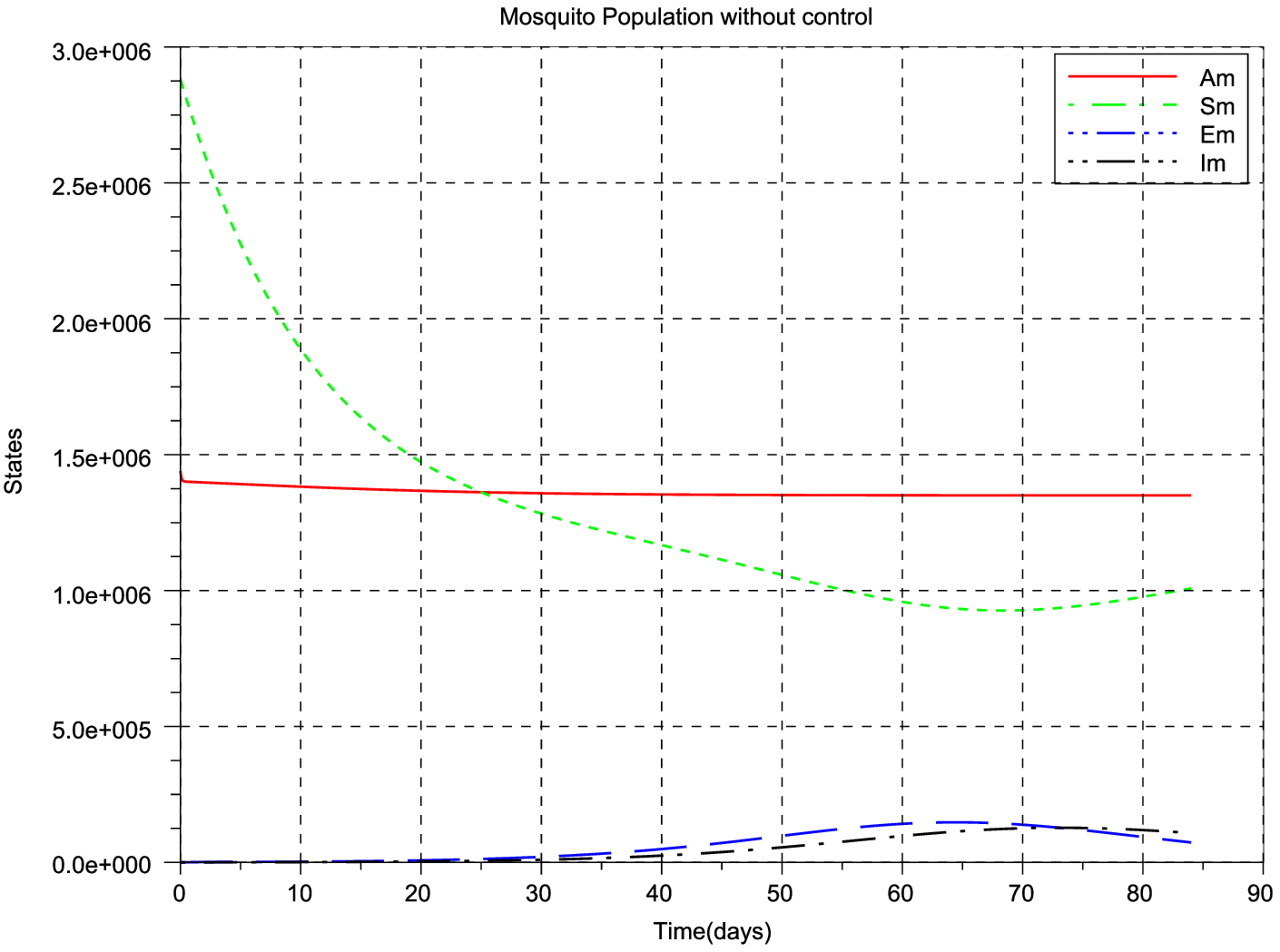}
{\caption{\label{mosquito_nocontrol} \small Mosquito compartments with no control}}
\end{minipage}
\end{figure}

The Figures~\ref{mosquito_control} and \ref{mosquito_nocontrol} 
show the difference between a region with control and without control.


The number of infected mosquitoes is close to zero in a situation 
where control is present. Note that we do not intend to eradicate 
the mosquitoes but instead the number of infected mosquitoes.


\section{Conclusions}
\label{sec:5}

It is very difficult to control or eliminate the \emph{Aedes aegypti} 
mosquito because it makes adaptations to the environment and becomes 
resistant to natural phenomena (\textrm{e.g.} droughts) 
or human interventions (\textrm{e.g.} control measures).

During outbreaks emergency vector control measures can also include 
broad application of insecticides. It has been shown here that with 
a steady spray campaign it is possible to reduce the number 
of infected humans and mosquitoes. Active monitoring and surveillance 
of the natural mosquito population should accompany control efforts 
to determine programme effectiveness.


\section*{Acknowledgements}

Work partially supported by Portuguese Foundation for Science 
and Technology (FCT) through the PhD Grant SFRH/BD/33384/2008 
(Rodrigues) and the R\&D units Algoritmi (Monteiro) and CIDMA (Torres).



\end{document}